\documentclass[10pt,a4paper,english]{article}
\usepackage[latin1]{inputenc}
\usepackage{babel}
\usepackage{amsmath,amsfonts,amssymb,amsthm}
\usepackage{bbm}
\usepackage{hyperref}
\usepackage{verbatim}
\usepackage{tikz}
\usepackage{todonotes}

\usetikzlibrary{positioning}

\theoremstyle{plain}
\newtheorem{teo}{Theorem}[section]
\newtheorem{cor}[teo]{Corollary}
\newtheorem{prop}[teo]{Proposition}
\newtheorem{lem}[teo]{Lemma}

\theoremstyle{defin}
\newtheorem{defin}[teo]{Definition}

\newcommand{\system}[1]{\mbox{\fontfamily{cmss}\fontshape{n}\fontseries{m}\selectfont#1}}
\newcommand{\ZF}{\system{ZF}}
\newcommand{\ZFC}{\system{ZFC}}
\newcommand{\AC}{\system{AC}}

\newcommand{\DC}{\system{DC}}

\newcommand{\pre}[2]{{}^{#1}#2}

\newcommand{\subsetsim}{\mathrel{\ooalign{\raise.4ex\hbox{$\subset$}\cr$\raise-.9ex\hbox{$\sim$}$}}}

\DeclareMathOperator{\cof}{cof}

\DeclareMathOperator{\Ord}{Ord}

\DeclareMathOperator{\ran}{ran}

\DeclareMathOperator{\lh}{lh}

\title{A Solovay-like model for singular generalized descriptive set theory}
\author{Vincenzo Dimonte\footnote{Universit\`{a} degli Studi di Udine, via delle Scienze, 206 33100 Udine (UD) \emph{E-mail address:} \texttt{vincenzo.dimonte@uniud.it}}}

\begin{document}

\maketitle

\begin{abstract}
 Kunen's proof of the non-existence of Reinhardt cardinals opened up the research on very large cardinals, i.e., hypotheses at the limit of inconsistency. One of these large cardinals, I0, proved to have descriptive-set-theoretical characteristics, similar to those implied by the Axiom of Determinacy: if $\lambda$ witnesses I0, then there is a topology for $V_{\lambda+1}$ that is completely metrizable and with weight $\lambda$ (i.e., it is a $\lambda$-Polish space), and it turns out that all the subsets of  $V_{\lambda+1}$ in $L(V_{\lambda+1})$ have the $\lambda$-Perfect Set Property in such topology. In this paper, we find another generalized Polish space of singular weight $\kappa$ of cofinality $\omega$ such that all its subsets have the $\kappa$-Perfect Set Property, and in doing this, we are lowering the consistency strength of such property from I0 to $\kappa$ $\theta$-supercompact, with $\theta>\kappa$ inaccessible.
\end{abstract}

What is Kunen's Theorem? Researchers in different areas will probably have different answers. In the field of very large cardinals, ``Kunen's Theorem'' is the fact that there are no elementary embeddings from the universe to itself (if we assume the Axiom of Choice) \cite{Kunen1971}. This was a surprising result, the first really non-obvious proof that a large cardinal is inconsistent, but the consequences of such theorem were equally unexpected.

Kunen's Theorem actually proves that there are no elementary embeddings from $V_{\lambda+2}$ to itself, for any $\lambda$. Looking for inconsistencies, a new breed of large cardinals was introduced, the rank-into-rank embeddings, starting from the existence of an elementary embedding from $V_{\lambda}$ to itself. Woodin pushed this further, defining I0, i.e., the existence of an elementary embedding from $L(V_{\lambda+1})$ to itself, with critical point less than $\lambda$. Instead of producing inconsistencies, I0 proved to be fruitful: it was initially used to prove the consistency of AD, and the more Woodin worked on it, the more similarities with $\mathbb{R}$ emerged, some of them with a descriptive set-theoretical ``flavor'', especially mirroring the structure of $L(\mathbb{R})$ under AD (more on this in \cite{Woodin2011}, \cite{Cramer2017}, \cite{Dimonte2018}).

Taking cue from these descriptive aspects, with Luca Motto Ros we developed a generalized descriptive set theory on cardinals with countable cofinality (since if I0$(\lambda)$ holds, then $\cof(\lambda)=\omega$): A $\lambda$-Polish space is a completely metrizable topology with weight $\leq\lambda$, and $\lambda$-Polish spaces share with Polish spaces many properties and characteristics. 

In this paper, we are interested in sets that have a tree structure. Cramer (\cite{Cramer2015}) proved that, under I0$(\lambda)$, all the subsets of $V_{\lambda+1}$ in $L(V_{\lambda+1})$ have an analogue of the Perfect Set Property, the $\lambda$-Perfect Set Property. This is analogous to the fact that, if $L(\mathbb{R})\vDash AD$, then all the subsets of $\mathbb{R}$ in $L(\mathbb{R})$ have the PSP. Even more, in \cite{DimonteMottoRos} it is described how this is a consequence of the fact that, under I0$(\lambda)$, all the sets have a property that is the generalization of being weakly homogeneously Souslin.

Looking at the classical case, a doubt emerge: while it is true that if $L(\mathbb{R})\vDash AD$, then all the subsets of $\mathbb{R}$ in $L(\mathbb{R})$ have the PSP, this is not the best scenario, consistency-wise. Solovay proved that if there is an inaccessible cardinal, and $V[G]$ is the generic extension via the Levy collapse of the inaccessible to $\omega$, then in $HOD_{\mathbb{R}}^{V[G]}$ all the subsets of $\mathbb{R}$ have the PSP, and an inaccessible cardinal is of consistency strength much lower than $AD^{L(\mathbb{R})}$. This begs the question: can we do the same for the $\lambda$-PSP? Is I0$(\lambda)$ necessary for having an example of a $\lambda$-Polish space where all its subsets have the $\lambda$-PSP, or we can lower the consistency strength of the large cardinals involved?

This paper provides a positive answer for the second question. Starting with a $\theta$-supercompact cardinal $\kappa$, with $\theta$ inaccessible, if $V[G]$ is a generic extension via supercompact Prikry forcing, then $V[G]$ contains an inner model in which all the subsets of $\pre{\omega}{\kappa}$, the ``$\kappa$-Baire space'', have the $\kappa$-PSP, and we can extend this result to any $\kappa$-Polish space in the inner model. We are using the construction made by Kafkoulis in \cite{Kafkoulis1994} of a generalized Solovay model: the construction is fairly similar to Solovay's construction, but where Solovay uses the homogeneity of the collapse and its simple factorization, we need to investigate in detail the Prikry-like properties of the supercompact Prikry forcing and its quotients.

It is still an open problem the exact lower bound for the existence of a $\lambda$-Polish space, with $\lambda>\omega$, with all the subsets with the $\lambda$-PSP, but surely it is a large cardinal.\footnote{Thanks to Gabriel Goldberg and Ralf Schindler for pointing out that it should be at least $\omega$-many Woodin cardinals: with the techniques in \cite{SchindlerSteel2014} it should be possible to prove, more in general, that if $L(V_{\lambda+1})\nvDash AC$, then there is an inner model with $\omega$-many Woodin cardinals, and it is easy to construct, in ZFC, a set that has not the $\lambda$-PSP.}

The research for this paper was partially supported by the Italian PRIN 2017 Grant ``Mathematical Logic: models, sets, computability''.

\section{Preliminaries and notation}

\textbf{Forcing constructions} In this subsection we collect some results about forcing constructions, with the further intent of fixing the notation. For more details, see for example the introductions of \cite{Shelah1998} or \cite{Abraham2010}.

In the following, all forcing posets are separative.

If $\mathbb{P}$ is a forcing poset, and $\mathbb{Q}\subseteq\mathbb{P}$ is a dense subset of $\mathbb{P}$, then $(\mathbb{Q},\leq)$ is also a forcing poset and:
\begin{itemize}
 \item If $G$ is $\mathbb{P}$-generic, then $G\cap\mathbb{Q}$ is $\mathbb{Q}$-generic and $V[G]=V[G\cap\mathbb{Q}]$;
 \item If $H$ is $\mathbb{Q}$-generic, then $H^{\uparrow}=\{p\in\mathbb{P} : \exists q\leq p\ q\in H\}$ is $\mathbb{P}$-generic and $V[H]=V[H^{\uparrow}]$. 
 \item If $\tau$ is a $\mathbb{Q}$-name, then it is also a $\mathbb{P}$-name, and if $\tau$ is a $\mathbb{P}$-name, then there is a $\mathbb{Q}$-name $\sigma$ such that $\Vdash_{\mathbb{P}}\tau=\sigma$;
 \item If $p\in\mathbb{Q}$ and $\tau$ is a $\mathbb{Q}$-name, then $p\Vdash_{\mathbb{P}}\varphi(\tau)$ iff $p\Vdash_{\mathbb{Q}}\varphi(\tau)$. 
\end{itemize}

If $\mathbb{P}$, $\mathbb{Q}$ are forcing posets, then $\pi:\mathbb{P}\to\mathbb{Q}$ is a dense homomorphism iff for any $p,q\in\mathbb{P}$, if $p\leq q$ then $\pi(p)\leq\pi(q)$, if $p$ is incompatible with $q$ then $\pi(p)$ is incompatible with $\pi(q)$, and $\ran(\pi)$ is dense in $\mathbb{Q}$. Then:
\begin{itemize}
 \item If $G$ is $\mathbb{P}$-generic, then $H=\{q\in\mathbb{Q}:\exists p\in G\ \pi(p)\leq q\}$ is $\mathbb{Q}$-generic and $V[G]=V[H]$;
 \item If $H$ is $\mathbb{Q}$-generic, then $G=\{p\in\mathbb{P}:\pi(p)\in H\}$ is $\mathbb{P}$-generic and $V[G]=V[H]$;
 \item $\pi$ extends to $\mathbb{P}$-names: if $\tau$ is a $\mathbb{P}$-name, define by induction on the rank $\pi(\tau)=\{(\pi(\sigma),\pi(p)):(\sigma,p)\in\tau\}$;
 \item For any $p\in\mathbb{P}$, $\tau$ $\mathbb{P}$-name, $p\Vdash_{\mathbb{P}}\varphi(\tau)$ iff $\pi(p)\Vdash_{\mathbb{Q}}\varphi(\pi(\tau))$.
\end{itemize}

Given $\mathbb{P}$, $\mathbb{R}$ forcing posets we say that $\mathbb{R}$ projects into $\mathbb{P}$, $\mathbb{P}\triangleleft\mathbb{R}$, iff there is a $\pi:\mathbb{R}\to\mathbb{P}$ such that:
 \begin{itemize}
  \item $\pi$ is order-preserving and $\pi(1_{\mathbb{R}})=1_{\mathbb{P}}$;
	\item For every $q\in\mathbb{R}$ and $p'\in\mathbb{P}$ such that $p'\leq \pi(q)$ there is a $q'\leq q$ in $\mathbb{R}$ such that $\pi(q')=p'$.
 \end{itemize}
In such a case, if $g$ is $\mathbb{P}$-generic, we can define $\mathbb{R}/g=\{r\in\mathbb{R}:\pi(r)\in g\}$. Let $\dot{\mathbb{Q}}$ be the $\mathbb{P}$-name for $\mathbb{R}/g$. Then the map $i:\mathbb{R}\to\mathbb{P}*\dot{\mathbb{Q}}$, defined as $i(r)=(\pi(r),\check{r})$, where $\check{r}$ is the canonical $\mathbb{P}$-name for $r$, is a dense homomorphism. We have the following connections between the two forcing notions:
\begin{itemize}
 \item If $H$ is a $\mathbb{R}$-generic filter over $V$, then $\pi''H$ is a $\mathbb{P}$-generic filter over $V$, and $H$ is a $\mathbb{R}/\pi''H$-generic filter over $V[\pi''H]$; moreover, $V[H]=V[\pi''H][H]$;
 \item If $g$ is $\mathbb{P}$-generic over $V$ and $h$ is $\dot{\mathbb{Q}}_g$-generic over $V[g]$, then $\{r\in\mathbb{R}:i(r)\in g*h\}=\{r\in\mathbb{R}:\pi(r)\in g\wedge r\in h\}=h$ is $\mathbb{R}$-generic over $V$ and $V[g][h]=V[h]$;
 \item If $\tau$ is a $\mathbb{R}$-name, then $i(\tau)$ is a $\mathbb{P}*\dot{\mathbb{Q}}$-name, and for any $g$ $\mathbb{P}$-generic over $V$ and $h$ $\dot{\mathbb{Q}}_g$-generic over $V[g]$ $(i(\tau)_g)_h=\tau_h$; 
 \item If $g$ is $\mathbb{P}$-generic, $p\in\dot{\mathbb{Q}}_g$ and $\tau$ is a $\dot{\mathbb{Q}}_g$-name, then $p\Vdash_{\dot{\mathbb{Q}}_g}\varphi(\tau)$ iff $p\Vdash_{\mathbb{R}}\varphi(\tau)$;
 \item For any $p\in\mathbb{R}$, $p\Vdash_{\mathbb{R}}\varphi(\tau)$ iff $\pi(p)\Vdash_{\mathbb{P}}(\check{p}\Vdash_{\dot{\mathbb{Q}}}\varphi(i(\tau)))$;
 \item For any $p\in\mathbb{P}$, $q\in\mathbb{R}$, $p\Vdash_{\mathbb{P}}(\check{q}\in\dot{\mathbb{Q}})$ iff $p\leq\pi(q)$.
\end{itemize}

\textbf{$\kappa$-Polish spaces}. In this subsection we collect some definitions and properties about $\kappa$-Polish spaces, as introduced in \cite{DimonteMottoRos}, where the base theory is ZF+$\AC_\kappa$. Since we will work in a model of ZF+$\DC_\kappa$, everything in this section will apply.

A topological space $X$ is $\kappa$-Polish if it is completely metrizable and has weight at most $\kappa$. If $\cof(\kappa)=\omega$, then $\pre{\omega}{\kappa}$, as the product of the discrete topologies on $\kappa$, is $\kappa$-Polish. If $(\eta_n : n\in\omega)$ is cofinal in $\kappa$, then also $C(\kappa)=\prod_{n\in\omega}\eta_n$, as the product of the discrete topologies on $\eta_n$, is $\kappa$-Polish. In fact, $\pre{\omega}{\kappa}\cong C(\kappa)$. The space $\pre{\kappa}{2}$, with the bounded topology (i.e., the topology generated by the basic open sets $N_s=\{x\in\pre{\kappa}{2}:x\upharpoonright\lh(s)=s\}$, with $s\in\pre{<\kappa}{2}$) is a completeley metrizable space with weight $2^{<\kappa}$ therefore if $2^{<\kappa}=\kappa$, then $\pre{\kappa}{2}$ is $\kappa$-Polish and isomorphic to $\pre{\omega}{\kappa}$.

If $\cof(\kappa)=\omega$, then the Woodin topology on $V_{\kappa+1}$ is the topology generated by the basic open sets $N_{(\alpha,a)}=\{A\subseteq V_\kappa : A\cap V_\alpha=a\}$, with $\alpha<\kappa$ and $a\subseteq V_\alpha$. If $|V_\kappa|=\kappa$, i.e., if $\kappa$ is a fixed point of the beth-function (for example if $\kappa$ is limit of inaccessible cardinals), then the space $V_{\kappa+1}$ with the Woodin topology is a $\kappa$-Polish space, and it is isomorphic to $\pre{\omega}{\kappa}$.

Let $X$ be a $\kappa$-Polish space, with $\cof(\kappa)=\omega$. Then there is a closed set $F\subseteq\pre{\omega}{\kappa}$ and a continuous bijection $f:F\to X$

If $X$ is a $\kappa$-Polish space, we say that $A\subseteq X$ has the $\kappa$-Perfect Set Property, or $\kappa$-PSP, if either $|A|\leq\kappa$, or $\pre{\kappa}{2}$ embeds into $A$ as a closed-in-$X$ set. If $2^{<\kappa}=\kappa$, this is equivalent to $\pre{\omega}{\kappa}$, or $C(\kappa)$, embedding into $A$ as a closed-in-$X$ set. In fact, this is equivalent to $\pre{\kappa}{2}$, $\pre{\omega}{\kappa}$, or $C(\kappa)$ embedding into $A$, disregarding the closure.

\section{Supercompact Prikry forcing, basic notions}

\begin{defin}
 A cardinal $\kappa$ is $\delta$-supercompact, with $\delta\geq\kappa$, iff there exists a normal, fine ultrafilter $U_\delta$ over ${\cal P}_\kappa(\delta)$. Such ultrafilter is called a $\delta$-supercompactness measure for $\kappa$.

 For $P,Q\in{\cal P}_\kappa(\delta)$, we say that $P\subsetsim Q$ (strong inclusion) iff $P\subset Q$ and $|P|<|Q\cap\kappa|$.

 For $B\subseteq{\cal P}_\kappa(\delta)$, we denote by $[B]^{[n]}$ the set of all $\subsetsim$-increasing $n$-length sequences of elements of $B$, by $[B]^{[<\omega]}$ the set of all $\subsetsim$-increasing finite sequences of elements of $B$, and by $[B]^{[\omega]}$ the set of all $\subsetsim$-increasing $\omega$-length sequences of elements of $B$.

 If $A\subseteq{\cal P}_\kappa(\delta)$ and $s\in[{\cal P}_\kappa(\delta)]^{[n]}$, then we call $A\setminus s=\{P\in A : s(n-1)\subsetsim P\}$. If $A\in U_\delta$, then $A\setminus s\in U_\delta$.

 If $s,t\in[{\cal P}_\kappa(\delta)]^{[<\omega]}$, we say $s\subsetsim t$ when $s(\lh(s)-1)\subsetsim t(0)$.
\end{defin}

\begin{prop}
 \label{prop:URowbottom}
 Let $\kappa$ be $\delta$-supercompact, and let $U_\delta$ be a normal, fine ultrafilter over ${\cal P}_\kappa(\delta)$. Let $\{A_P : P\in{\cal P}_\kappa(\delta)\}$ be a family of sets in $U_\delta$. Then 
\begin{equation*}
 \triangle_{P\in{\cal P}_\kappa(\delta)}A_P=\{Q\in{\cal P}_\kappa(\delta) : \forall P\subsetsim Q\ Q\in A_P\}\in U_\delta.
\end{equation*}
 Moreover, if $F:[{\cal P}_\kappa(\delta)]^{[<\omega]}\to 2$, then there is an $A\in U_\delta$ such that for each $n\in\omega$, $F$ is constant on $[A]^{[n]}$. 
\end{prop}

\begin{defin}
 Let $\kappa$ be a $\delta$-supercompact cardinal, and let $U_\delta$ be a $\delta$-supercompact measure for $\kappa$. 

 The elements of $\mathbb{P}_{U_\delta}$, the supercompact Prikry forcing for $U_\delta$, are all the sets of the form $p=(s,A)$, with $s=(P_0,\dots,P_{n-1})$ a finite $\subsetsim$-increasing sequence of elements of ${\cal P}_\kappa(\delta)$, and $A\in U_\delta$ such that for all $Q\in A$, $P_{n-1}\subsetsim Q$. The sequence $s$ is called the stem of $p$, and we denote it as $stem(p)$. We call $\lh(p)=\lh(stem(p))$. The set $A$ is called the measure-one part of $p$.

 We say that $(t,B)\leq(s,A)$, with $(s,A),(t,B)\in\mathbb{P}_{U_\delta}$, iff $s\sqsubseteq t$ (i.e., $\lh(t)\geq\lh(s)$ and for all $i<\lh(s)$ $t(i)=s(i)$), for all $\lh(s)\leq i<\lh(t)$ $t(i)\in A$, and $B\subseteq A$. 

 We say that $(t,B)\leq^*(s,A)$ iff $(t,B)\leq(s,A)$ and $s=t$.
\end{defin}

Let $\mathbb{P}_{U_\delta}$ be the supercompact Prikry forcing for $U_\delta$. Every $\mathbb{P}_{U_\delta}$-generic $G$ induces an infinite $\subsetsim$-increasing sequence $\sigma_G=\bigcup\{stem(p) : p\in G\}$, that by density arguments is such that for each $\alpha\leq\delta$, $\alpha=\bigcup_{n\in\omega}(\sigma_G(n)\cap\alpha)$. 

Vice versa, if we have a $\sigma\in[{\cal P}_\kappa(\delta)]^{[\omega]}$, we can define the filter in $\mathbb{P}_{U_\delta}$ generated by $\sigma$: 
\begin{equation*}
 {\cal F}_{\sigma}=\{(t,A)\in\mathbb{P}_{U_\delta} : t=\sigma\upharpoonright\lh(t),\ \forall i\geq\lh(t)\ \sigma(i)\in A\}.
\end{equation*}

It is a maximal filter, but it is not necessarily generic.

\begin{lem}
 \label{lem:filterfromsequence}
 Let $\mathbb{P}_{U_\delta}$ be the supercompact Prikry forcing for $U_\delta$ and let $G$ be $\mathbb{P}_{U_\delta}$-generic. Then ${\cal F}_{\sigma_G}=G$.
\end{lem}
\begin{proof}
 Let $p=(s,A)\in G$, $\lh(p)=n$. First of all, for each $i<n$ $\sigma_G(i)=s(i)$, therefore $s=\sigma_G\upharpoonright n$. Since for every $m\geq n$, $D_m=\{p\in\mathbb{P}_{U_\delta} : \lh(p)>m\}$ is dense, then for each $m\geq n$ there is a $p_m\leq p$, $p_m\in G$ such that $\lh(p_m)>m$. So if $n\leq i<m$ then $p_m(i)=\sigma_G(i)\in A$, i.e., $p\in{\cal F}_{\sigma_G}$.

We proved that $G\subseteq{\cal F}_{\sigma_G}$. But since $G$ is maximal, $G={\cal F}_{\sigma_G}$.
\end{proof}

\begin{teo}
 \label{thm:basic}
 Let $\kappa$ be $\delta$-supercompact, with $\cof(\delta)>\kappa$, and let $U_\delta$ be a $\delta$-supercompactness measure for $\kappa$. Let $G$ be a generic set for $\mathbb{P}_{U_\delta}$ over $V$. Then in $V[G]$ no new bounded subsets are added to $\kappa$, every cardinal in the interval $[\kappa,\delta]$ of cofinality $\geq\kappa$ in $V$ has cofinality $\omega$, and all the cardinals above $\delta$ are preserved. Therefore, in $V[G]$, $|\delta|=\kappa$.
\end{teo}

Because of Proposition \ref{prop:URowbottom}, supercompact Prikry forcing enjoys the usual Prikry-like properties:

\begin{prop}[Prikry condition]
 Let $\kappa$ be $\delta$-supercompact, $U_\delta$ be a $\delta$-supercompactness measure for $\kappa$ and $\mathbb{P}_{U_\delta}$ be the supercompact Prikry forcing for $U_\delta$. Let $\psi$ be a forcing statement and $p\in\mathbb{P}_{U_\delta}$. Then there exists $p'\leq^* p$ that decides $\psi$, i.e., either $p'\Vdash\psi$ or $p'\Vdash\neg\psi$.
\end{prop}

\begin{prop}[Geometric condition]
  \label{prop:geometric}
  For all $\sigma\in[{\cal P}_\kappa(\delta)]^{[\omega]}$, the filter ${\cal F}_\sigma$ is generic iff for all $A\in U_\delta$ there is an $n\in\omega$ such that $\sigma\upharpoonright [n,\omega)\subseteq A$. 
\end{prop}

\section{Supercompact Prikry forcing, quotients}

In this section, $\kappa$ will always be a $\delta$-supercompact cardinal with $\cof(\delta)>\kappa$, $U_\delta$ will be a $\delta$-supercompactness measure, and $\alpha\in[\kappa,\delta)$.

\begin{defin}
 Let $s\in[{\cal P}_\kappa(\delta)]^{[n]}$, and $\alpha<\delta$. Then $s\downarrow\alpha=(s(0)\cap\alpha,\dots,s(n-1)\cap\alpha)$. Clearly, $s\downarrow\alpha\in[{\cal P}_\kappa(\alpha)]^{[n]}$.

 Given a $\delta$-supercompactness measure $U_\delta$ for $\kappa$, if $A\in U_\delta$, then let $A\downarrow\alpha=\{P\cap\alpha : P\in A\}$. Define also $U_\alpha=\{A\downarrow\alpha : A\in U_\delta\}$. It is an $\alpha$-supercompactness measure for $\kappa$. Moreover, if $\alpha<\beta<\delta$ then $U_\beta\downarrow\alpha=U_\alpha$.
\end{defin}

Note that if $B\in U_\alpha$, then $\{P\in{\cal P}_\kappa(\delta) : P\cap\alpha\in B\}\in U_\delta$, i.e., $U_\delta$ projects into $U_\alpha$.

\begin{defin}
 Let $\kappa$ be $\delta$-supercompact, $U_\delta$ be a $\delta$-supercompactness measure for $\kappa$ and $\mathbb{P}_{U_\delta}$ be the supercompact Prikry forcing for $U_\delta$. 

 If $p=(s,A)\in \mathbb{P}_{U_\delta}$, then $p\downarrow\alpha=(s\downarrow\alpha,A\downarrow\alpha)\in\mathbb{P}_{U_\alpha}$.

 If $G$ is $\mathbb{P}_{U_\delta}$-generic, then $G\downarrow\alpha=\{p\downarrow\alpha : p\in G\}$. 
\end{defin}

When $\alpha<\delta$, we would like to see $\mathbb{P}_{U_\alpha}$ as a projection of $\mathbb{P}_{U_\delta}$, but this is not technically true. In fact, to have $\mathbb{P}_{U_\alpha}\vartriangleleft\mathbb{P}_{U_\delta}$, we should have that for each $q\in\mathbb{P}_{U_\delta}$ and for each $p\in\mathbb{P}_{U_\alpha}$ such that $p\leq q\downarrow\alpha$ there exists a $q'\leq q$ such that $q'\downarrow\alpha=p$, and this is not true for all the conditions of $\mathbb{P}_{U_\delta}$. The key obstacle is that if $t\in[A\downarrow\alpha]^{[<\omega]}$, with $A\in U_\delta$, there is an $r\in [A]^\omega$ such that $r\downarrow\alpha$, but it is not necessary that $r$ is $\subsetsim$-increasing. But if we fix $\alpha$, we can find a dense subforcing for which this is true.

\begin{defin}
  Let $\alpha\in[\kappa,\delta)$. We say that $p=(s,A)\in\mathbb{P}_{U_\delta}$ is $\alpha$-nice iff for all $t\in[A\downarrow\alpha]^{[<\omega]}$ there is an $r\in[A]^{[<\omega]}$ such that $r\downarrow\alpha=t$.
\end{defin}

\begin{lem}{\cite[Lemma 2.1.9]{Kafkoulis1994}}
 \label{lem:nicedense}
 Let $\alpha\in[\kappa,\delta)$. Then for any $p\in\mathbb{P}_{U_\delta}$ there is a $p'\leq^* p$ that is $\alpha$-nice. Therefore the set of $\alpha$-nice conditions is dense in $\mathbb{P}_{U_\delta}$.
\end{lem}

Call $\mathbb{P}_{U_\delta}^\alpha=\{p\in\mathbb{P}_{U_\delta} : p$ is $\alpha$-nice$\}$. If $G$ is $\mathbb{P}_{U_\delta}$-generic, then $G\cap\mathbb{P}_{U_\delta}^\alpha$ is $\mathbb{P}_{U_\delta}^\alpha$-generic, and $V[G]=V[G\cap\mathbb{P}_{U_\delta}^\alpha]$. On the other hand, if $H$ is $\mathbb{P}_{U_\delta}^\alpha$-generic, then $H^{\uparrow}=\{p\in\mathbb{P}_{U_\delta}:\exists q\leq p\ q\in H\}$ is $\mathbb{P}_{U_\delta}$-generic and $V[H]=V[H^{\uparrow}]$. Also, for any $\mathbb{P}_{U_\delta}$-name $\tau$ there is a $\mathbb{P}_{U_\delta}^\alpha$-name $\tau^\alpha$ such that $\Vdash_{\mathbb{P}_{U_\delta}}\tau=\tau^\alpha$.  

We defined $\mathbb{P}_{U_\delta}^\alpha$ so that $\mathbb{P}_{U_\alpha}\vartriangleleft\mathbb{P}_{U_\delta}^\alpha$, as witnessed by $\downarrow\alpha$. The first consequence is that if $G$ is a $\mathbb{P}_{U_\delta}$-generic filter, then $H=G\cap\mathbb{P}_{U_\delta}^\alpha$ is a $\mathbb{P}_{U_\delta}^\alpha$-generic filter, and $g=H\downarrow\alpha$ is $\mathbb{P}_{U_\alpha}$-generic. 

Then we can quotient $\mathbb{P}_{U_\delta}^\alpha$ via $\mathbb{P}_{U_\alpha}$: if $g$ is $\mathbb{P}_{U_\alpha}$-generic, let $\mathbb{Q}_{\alpha\delta}(g)=\{q\in \mathbb{P}_{U_\delta}^\alpha : q\downarrow\alpha\in g\}=\{q\in\mathbb{P}_{U_\alpha} : q$ is $\alpha$-nice and $q\downarrow\alpha\in g\}=\mathbb{P}_{U_\delta}^\alpha/g$. Let $\dot{\mathbb{Q}}_{\alpha\delta}$ be the $\mathbb{P}_{U_\alpha}$-name for $\mathbb{Q}_{\alpha\delta}(g)$. Then the map $i:\mathbb{P}_{U_\delta}^\alpha\to\mathbb{P}_{U_\alpha}*\dot{\mathbb{Q}}_{\alpha\delta}$, defined as $i(r)=(r\downarrow\alpha,\check{r})$, is a dense homomorphism. So if $H$ is a $\mathbb{P}_{U_\delta}^\alpha$-generic filter, then it is also $\mathbb{Q}_{\alpha\delta}(H\downarrow\alpha)$-generic over $V[H\downarrow\alpha]$, and $V[H]=V[H\downarrow\alpha][H]$. On the other hand, if $g$ is $\mathbb{P}_{U_\alpha}$-generic and $h$ is $\mathbb{Q}_{\alpha\delta}(g)$-generic over $V[g]$, then $h$ is also $\mathbb{P}_{U_\delta}^\alpha$-generic over $V$ and $V[g][h]=V[h]$.

With the above notation, the following is standard:

\begin{lem}
 \label{lem:correspondence}
  \begin{enumerate}
   \item Let $\tau$ be a $\mathbb{P}_{U_\delta}$-name. For each $p\in\mathbb{P}_{U_\delta}^\alpha$ the following are equivalent:
    \begin{itemize}
     \item $p\Vdash_{\mathbb{P}_{U_\delta}}\varphi(\tau)$;
	   \item $p\Vdash_{\mathbb{P}_{U_\delta}^\alpha}\varphi(\tau^\alpha)$
	   \item $p\downarrow\alpha\Vdash_{\mathbb{P}_{U_\alpha}}(\check{p}\Vdash_{\dot{\mathbb{Q}}_{\alpha\delta}}\varphi(i(\tau^\alpha)))$.
    \end{itemize}
   \item If $p\in\mathbb{Q}_{\alpha\delta}(g)$ for some $g$ $\mathbb{P}_{U_\alpha}$-generic and $\tau$ is a $\mathbb{Q}_{\alpha\delta}(g)$-name, then we can add
    \begin{itemize}
     \item $p\Vdash_{\mathbb{Q}_{\alpha\delta}(g)}\varphi(\tau)$.
    \end{itemize}
	\end{enumerate}
\end{lem}

\begin{prop}
 \label{prop:manyextensions}
 Let $g$ be $\mathbb{P}_{U_\alpha}$-generic over $V$, and let $\tau$ be a $\mathbb{Q}_{\alpha\delta}(g)$-name for a function in $\pre{\omega}{\kappa}$. Let $p\in\mathbb{Q}_{\alpha\delta}(g)$ be such that $p\Vdash\tau:\check{\omega}\to\check{\kappa}$. Consider $A_p=\{\gamma<\kappa : \exists q\leq p\ \exists m\in\omega\ q\Vdash\tau(\check{m})=\check{\gamma}\}$. If $|A_p|<\kappa$ in $V[g]$, then for each $h$ $\mathbb{Q}_{\alpha\delta}(g)$-generic such that $p\in h$, $\tau_h\in V[g]$. \footnote{We thank the referee for having greatly improved the proof of this proposition.}
\end{prop}
\begin{proof}
 (In this proof we are only using the forcing relation $\Vdash_{\mathbb{Q}_{\alpha\delta}(g)}$, so we just write $\Vdash$).

 Since $|A_p|<\kappa$ in $V[g]$, then also $B_p=\{(m,\gamma):\exists q\leq p,\ q\in\mathbb{Q}_{\alpha\delta}(g)\wedge q\Vdash\tau(m)=\gamma\}$ has cardinality less than $\kappa$ in $V[g]$. Let $\eta<\kappa$ be such that $(|B_p|=\eta)^{V[g]}$, and let $f\in V[g]$ be a bijection that witnesses it. 

 Let $h$ be $\mathbb{Q}_{\alpha\delta}(g)$-generic such that $p\in h$. Then $\tau_h\subseteq B_p$, and $f''\tau_h\subseteq\eta$. But then $f''\tau_h$ is a bounded subset of $\kappa$ in $V[g][h]$, that is a $\mathbb{P}_{U_\delta}$-generic extension, therefore by Theorem \ref{thm:basic} $f''\tau_h\in V$. Since $f$ is a bijection, $\tau_h$ is definable from $f$ and $f''\tau_h$, so $\tau_h\in V[g]$.
\end{proof}

We are going now to approach the problem of homogeneity for the quotient forcing.

\begin{teo}{\cite[Before Lemma 3.4]{Magidor1977}}
 Each $f$ permutation of $\delta$  such that $f\upharpoonright\kappa= id$ induces an automorphism $\pi_f$ of $\mathbb{P}_{U_\delta}$ in the following way: If $s\in[{\cal P}_\kappa(\delta)]^{[<\omega]}$, then $\pi_f(s)=(f''s(0),\dots,f''s(\lh(s)-1))$, and if $p=(s,A)\in\mathbb{P}_{U_\delta}$, we define $\pi_f(p)=(\pi_f(s),\pi_f''A)$. 
\end{teo}

If $f$ is a permutation such that $f\upharpoonright\alpha=id$, then $\pi_f(p)\downarrow\alpha=p\downarrow\alpha$, and $f$ induces an automorphism of $\mathbb{Q}_{\alpha\delta}(g)$. This consideration permits us to prove that the forcing $\mathbb{Q}_{\alpha\delta}(g)$ is weakly homogeneous, as witnessed by the $\pi_f$'s:

\begin{lem}{\cite[Lemma 2.1.10]{Kafkoulis1994}}
 \label{lem:homogeneity}
 Let $\alpha<\delta$, $p\in\mathbb{P}^\alpha_{U_\delta}$, $q\in\mathbb{P}_{U_\delta}$ be such that $p\upharpoonright\alpha=q\upharpoonright\alpha$. Then there is a permutation $f$ of $\delta$ such that $f\upharpoonright\alpha=id$, and for all $q_1\leq q$ there is a $p_1\leq p$ such that $p_1\upharpoonright\alpha=q_1\upharpoonright\alpha$ and $\pi_f(p_1)$ is compatible with $q_1$. 
\end{lem}

\begin{cor}
 \label{cor:homogeneity}
 Let $\alpha<\delta$, $g$ $\mathbb{P}_{U_\alpha}$-generic over $V$, let ${\cal A}=\{\pi_f:f$ is a permutation of $\delta$ such that $f\upharpoonright\kappa=id\}$. Then $\mathbb{Q}_{\alpha\delta}(g)$ is weakly homogeneous, and the homogeneity is witnessed by elements in ${\cal A}$. In particular, if $\tau$ is a $\mathbb{Q}_{\alpha\delta}(g)$-name that is invariant under members of ${\cal A}$ (for example, a canonical check-name for an element of $V[g]$), then for any formula $\varphi$, if $q\Vdash_{\mathbb{Q}_{\alpha\delta}(g)}\varphi(\tau)$ then $\Vdash_{\mathbb{Q}_{\alpha\delta}(g)}\varphi(\tau)$.
\end{cor}
\begin{proof}
 Let $p,q\in\mathbb{Q}_{\alpha\delta}(g)$. Then $p\downarrow\alpha$ and $q\downarrow\alpha$ are in $g$, therefore there is an $r\in g$ such that $r\leq p\downarrow\alpha,\ q\downarrow\alpha$. Since $p,q\in\mathbb{P}^\alpha_{U_\delta}$ and $\mathbb{P}_{U_\alpha}\triangleleft\mathbb{P}^\alpha_{U_\delta}$, there are $p',q'\in\mathbb{P}^\alpha_{U_\delta}$ such that $p'\leq p$, $q'\leq q$ and $p'\downarrow\alpha=r=q'\downarrow\alpha$. Since $r\in g$, then, $p',q'\in\mathbb{Q}_{\alpha\delta}(g)$. By Lemma \ref{lem:homogeneity} there exists a permutation $f$ of $\delta$ such that $f\upharpoonright\alpha=id$ and a $p''\leq p'$ such that $\pi_f(p'')$ is compatible with $q'$. But then also $\pi_f(p)$ is compatible with $q$. 
\end{proof}

\section{The main construction}

In this section, $\kappa$ will always be a $\theta$-supercompact cardinal, with $\theta>\kappa$ an inaccessible cardinal. To avoid cluttering, we are skipping $\theta$ in the notations, therefore $U$ will be a $\theta$-supercompactness measure, and if $\alpha\in[\kappa,\theta)$, then $\dot{\mathbb{Q}}_{\alpha}$ is the quotient of $\mathbb{P}^\alpha_U$ over $\mathbb{P}_{U_\alpha}$.

Let $G$ be $\mathbb{P}_{U}$-generic over $V$. Then we define, in $V[G]$, 
\begin{equation*}
 H^G=\bigcup\{{\cal P}(\kappa)\cap V[G\downarrow\alpha]|\alpha<\theta\}. 
\end{equation*}

Our model of reference is going to be $L(H^G)$.

\begin{prop}{\cite[Lemma 2.1.11, Corollary 2.1.12, Corollary 2.1.13]{Kafkoulis1994}}
 \label{prop:basicH}
  In $V[G]$, $H^G={\cal P}(\kappa)\cap L(H^G)$, so for each $z\in{\cal P}(\kappa)\cap L(H^G)$ there is an $\alpha<\theta$ such that $z\in V[G\downarrow\alpha]$, and $L(H^G)\vDash\theta=\kappa^+$.
\end{prop}

\begin{lem}{\cite[Lemma 2.1.11]{Kafkoulis1994}}
 \label{lem:invariantname}
 Let $\beta<\theta$. In $V[G]$, let $\dot{g}_\beta=\{(\check{P},p):p=((P_1,\dots,P_n,A)\in\mathbb{P}_{U}\wedge\exists i\ P=P_i\cap\beta\}$. Then $\dot{g}_\beta$ is a $\mathbb{P}_{U}$-name for $G\downarrow\beta$. Let ${\cal A}$ be the set of automorphisms generated by permutations $f$ of $\theta$ such that $f\downarrow\kappa=id$. Then
 \begin{equation*}
  H^G=\bigcup\{{\cal P}(\kappa)\cap V[(h(\dot{g}_\beta))_G] : h\in{\cal A}\wedge\beta\in[\kappa,\theta)\}
 \end{equation*} 
 and therefore there is a $\mathbb{P}_{U}$-name for $H^G$ that is invariant under members of ${\cal A}$.
\end{lem}

Let $\alpha<\theta$. Then $|tr(\mathbb{P}_{U_\alpha})|\leq 2^{|\alpha|^+}$. Let $\eta=(2^{|\alpha|^+})^V$. Then, since $\theta$ is inaccessible, there is a $\beta<\theta$ such that $\beta>\eta$, and in $V[G\downarrow\beta]$ we have that $|tr(\mathbb{P}_{U_\alpha})|\leq\kappa$. But then it is possible to define, in $V[G\downarrow\beta]$, a set $E\subseteq\kappa\times\kappa$ such that $(\mathbb{P}_\alpha,\in)$ is isomorphic to $(\kappa,E)$. So $E$ is codeable as a member of ${\cal P}(\kappa)$, therefore it is in ${\cal P}(\kappa)\cap V[G\downarrow\beta]$ and so in $H^G$, therefore $\mathbb{P}_{U_\alpha}\in L(H^G)$. In fact, the same holds for any element of $(H_{\kappa^+})^{V[G\downarrow\alpha]}$ for any $\alpha<\theta$, so in particular, for all $\alpha<\theta$ then $U_\alpha,\mathbb{P}_{U_\alpha},G\downarrow\alpha,(\pre{\omega}{\kappa})^{V[G\downarrow\alpha]}\in L(H^G)$, and for any $\alpha<\beta$ $\mathbb{Q}_{\alpha\beta}(G\downarrow\alpha)\in L(H^G)$.

Notice also that, since $(\cof(\kappa)=\omega)^{V[G\downarrow\kappa]}$, then there is a $\omega$-sequence in ${\cal P}(\kappa)\cap V[G\downarrow\kappa]$ that is cofinal in $\kappa$, therefore such a sequence is in $H^G$, and so $(\cof(\kappa)=\omega)^{L(H^G)}$. 

\begin{teo}{\cite[Theorem 4.1.13]{Kafkoulis1994}}
 $L(H^G)\vDash\DC_\kappa$.
\end{teo}

Therefore we can develop in $L(H^G)$ a generalized descriptive set theory as in \cite{DimonteMottoRos}, $\pre{\omega}{\kappa}$ is a $\kappa$-Polish space in $L(H^G)$, and it makes sense to ask in $L(H^G)$ which subsets of $\pre{\omega}{\kappa}$ have the $\kappa$-PSP. Note also that since $\kappa$ is $\theta$-supercompact, $|V_\kappa|=\kappa$, and this is true also in $L(H^G)$, therefore $(V_{\kappa+1})^{L(H^G)}$ is a $\kappa$-Polish space in $L(H^G)$ homeomorphic to $(\pre{\omega}{\kappa})^{L(H^G)}$, and in fact $L(H^G)=(L(V_{\kappa+1}))^{L(H^G)}$.

Our first objective will be to find a way to identify and construct $\mathbb{Q}_{\alpha\beta}(G\downarrow\alpha)$-generics in $L(H^G)$, for $\kappa\leq\alpha<\beta<\theta$. The following is a standard fact of Prikry-like forcings, (see e.g. \cite[Theorem 2.2.6]{Kafkoulis1994}), but we redo the proof because we need a slightly more precise statement. The proof uses the same ideas of Mathias' proof of the geometric condition.

\begin{prop}
 \label{prop:sequenceforforcing}
 Let $\beta<\theta$. Then there is a $\subset$-descending sequence $\langle A_n^* : n\in\omega\rangle$ in $L(H^G)$ such that $A^*_n\in U_\beta$, and for each $\sigma\in [{\cal P}_\kappa(\beta)]^{[\omega]}$, ${\cal F}_\sigma$ is $\mathbb{P}_{U_\beta}$-generic over $V$ iff for all $n\in\omega$ there is an $m_n\in\omega$ such that for all $i>m_n$, $\sigma(i)\in A_n^*$. 
\end{prop}
\begin{proof}
 Since $\theta$ is inaccessible, there exists a $\gamma<\theta$ such that $(|U_\beta|=\kappa)^{V[G\downarrow\gamma]}$, and by usual considerations (see e.g. \cite[Corollary 2.16]{Prikry1970}) there exists in $V[G\downarrow\gamma]$ a collection of $E_n\in V$ such that $U_\beta=\bigcup_{n\in\omega}E_n$ and $|E_n|<\kappa$. Since $\langle E_n:n\in\omega\rangle\in (H_{\kappa^+})^{V[G\downarrow\gamma]}$, then $\langle E_n:n\in\omega\rangle\in L(H^G)$. Let $E^*_m=\bigcup_{n\leq m} E_k$ and $A^*_m=\bigcap E^*_m$. Since $E^*_m\in V$, by the $\kappa$-completeness of $U_\beta$ $A^*_m\in U_\beta$, and $\langle A^*_m:m\in\omega\rangle\in L(H^G)$. 

 Let $\sigma\in [{\cal P}_\kappa(\beta)]^{[\omega]}$. If ${\cal F}_\sigma$ is $\mathbb{P}_{U_\beta}$-generic over $V$, then by the Geometric Condition (Proposition \ref{prop:geometric}) for any $A\in U_\beta$ there is a $n_A\in\omega$ such that for all $i>n_A$, $\sigma(i)\in A$, and this is true also for the $A_n^*$'s. On the other hand, if for all $n\in\omega$ there is an $m_n\in\omega$ such that for all $i>m_n$, $\sigma(i)\in A_n^*$, let $A\in U_\beta$ and let $n\in\omega$ be such that $A\in E^*_n$. Then $A^*_n\subseteq A$, so there is an $m_n\in\omega$ such that for all $i>m_n$, $\sigma(i)\in A$, therefore by the Geometric Condition ${\cal F}_\sigma$ is $\mathbb{P}_{U_\beta}$-generic over $V$.
\end{proof}

Now we introduce the main technique, that assures us that we can extend any $\mathbb{P}_{U_\alpha}$-generic to a $\kappa$-perfect set of $\mathbb{P}_{U_\beta}$-generics.

\begin{teo}
 \label{thm:perfect}
 Let $\alpha<\beta<\theta$, let $g=G\downarrow\alpha$. Let $\tau$ be a $\mathbb{Q}_{\alpha\beta}(g)$-name for an element of $\pre{\omega}{\kappa}$ that is not in $V[g]$. Let $p\in\mathbb{Q}_{\alpha\beta}(g)$ be such that $p\Vdash\tau:\check{\omega}\to\check{\kappa}$. Then $\{\tau_h : h$ is $\mathbb{Q}_{\alpha\beta}(g)$-generic, $p\in h\}$, as calculated in $L(H^G)$, contains a $\kappa$-perfect set.  
\end{teo}
\begin{proof}
 (In this proof we use only the forcing relation $\Vdash_{\mathbb{Q}_{\alpha\beta}(g)}$, therefore to avoid cluttering we call it just $\Vdash$.)

 Since $g\in L(H^G)$, then $(\cof(\kappa)=\omega)^{L(H^G)}$. Fix $(\eta_n : n\in\omega)\in L(H^G)$, a sequence cofinal in $\kappa$. Fix also a sequence $(A^*_n : n\in\omega)\in L(H^G)$ for $\mathbb{P}_{U_\beta}$ as in Proposition \ref{prop:sequenceforforcing}. 

 The idea is the following: for each $s\in\bigcup_{n\in\omega}\prod_{i<n}\eta_i$, we are going to define in $L(H^G)$ $p_s,q_s,r_s\in\mathbb{Q}_{\alpha\beta}(g)$, so that $p_s\leq r_s\leq q_s$ and if $t=s^\smallfrown\langle\xi\rangle$, then $q_t\leq p_s$, in the following way:
 \begin{itemize}
  \item if $\lh(s)=n$, for each $\xi_1,\xi_2\in \eta_n$, if $\xi_1\neq\xi_2$ then $q_{s^\smallfrown\{\xi_1\}}$ and $q_{s^\smallfrown\{\xi_2\}}$ are going to force a different behavior for $\tau$, and for this we will use Proposition \ref{prop:manyextensions};
	\item if $x\in\prod_{i<\omega}\eta_i$, then $\langle p_{x\upharpoonright n}:n\in\omega\rangle$ generates a $\mathbb{Q}_{\alpha\beta}(g)$-generic filter, and for this we will use Proposition \ref{prop:sequenceforforcing};
	\item the condition $r_s$ decides the behavior of $\tau$ completely up to a certain point.
 \end{itemize}

 At the end of the construction, for each $x\in\prod_{i<\omega}\eta_i$, $\langle q_{x\upharpoonright n}:n\in\omega\rangle$ generates a $\mathbb{Q}_{\alpha\beta}(g)$-generic filter $G_x$ (because of the $p_{x\upharpoonright n}$'s). Then the function that to each $x\in\prod_{i<\omega}\eta_i$ associates the interpretation of $\tau$  under $G_x$ is injective (because of the way we have chosen the $q$'s) and continuous (because of the way we have chosen the $r$'s).

Let $p=p_\emptyset\in\mathbb{Q}_{\alpha\beta}(g)$. By Proposition \ref{prop:manyextensions} 
\begin{equation*}
 |\{\gamma<\kappa : \exists q\leq p_\emptyset\ \exists m\in\omega\ q\Vdash\tau(\check{m})=\check{\gamma}\}|=\kappa.
\end{equation*}
 in $V[g]$. Let 
\begin{equation*}
 \Gamma_m(p_\emptyset)=\{\gamma<\kappa : \exists q\leq p_\emptyset\ q\Vdash\tau(\check{m})=\check{\gamma}\}.
\end{equation*}

Then there must exist an $m_{\emptyset}$ such that $|\Gamma_{m_\emptyset}(p_\emptyset)|\geq\eta_0$. Let $\{\gamma_\xi : \xi<\eta_0\}$ be an enumeration of a subset of $\Gamma_{m_\emptyset}(p_\emptyset)$. Then for each $\xi<\eta_0$ choose $q_{(\xi)}\leq p_\emptyset$ such that $q_{(\xi)}\Vdash\tau(\check{m}_\emptyset)=\gamma_\xi$. The sequence $\langle q_{\langle\xi\rangle}:\xi<\eta_0\rangle$ has been constructed in $V[g]$, but by the remark after Lemma \ref{lem:invariantname} such sequence is also in $L(H^G)$.

Extend $q_{(\xi)}$ to $r_{(\xi)}$ so that for each $m<m_\emptyset$, $r_{(\xi)}$ decides the value of $\tau(m)$. 

So let $r_{(\xi)}=(s,A)$. Since $r_{(\xi)}\in\mathbb{Q}_{\alpha\beta}(g)$, $s\downarrow\alpha=\sigma_g\upharpoonright\lh(s)$, and $\sigma_g\upharpoonright[\lh(s),\omega)\subseteq A\downarrow\alpha$. We would like to use Proposition \ref{prop:sequenceforforcing} to find a $p=(t,B)\leq r_{(\xi)}$ such that $B\subseteq A_0^*$, but we cannot just take $(s, A\cap A_0^*)$, since it is possible that this condition is not in $\mathbb{Q}_{\alpha\beta}(g)$. We need to find in $L(H^G)$ a $\mathbb{Q}_{\alpha\beta}(g)$-generic set to use as a guide to define $p$.

Define $T$ in $L(H^G)$ in the following way:\footnote{The construction of $T$ is as in the proof of \cite[Lemma 2.2.7]{Kafkoulis1994}}

The tree $T$ is generated by the sequences $\langle P_0,\dots,P_m\rangle\in[{\cal P}_\kappa(\beta)]^{[<\omega]}$ such that $\langle P_0,\dots,P_m\rangle\upharpoonright\lh(s)=s$, for all $\lh(s)\leq i\leq n$ $P_i\in A$, there exist $n\in\omega$, $i_0,\dots,i_n$ such that if $j>i_k$ then $P_j\in A_k^*$, $P_m\in A_n^*$, and $\langle P_0,\dots,P_m\rangle\downarrow\alpha=\sigma_g\upharpoonright m$. Since the definition uses only the sequence $\langle A_n^*:n\in\omega\rangle$ and $g$, it is clear that $T\in L(H^G)$, and if $\sigma$ is a branch in $T$, then $s\sqsubseteq\sigma$, $\sigma\downarrow\alpha=\sigma_g$, for all $i\geq\lh(s)$ $\sigma(i)\in A$, and for all $n\in\omega$ there is an $m_n\in\omega$ such that for all $i\geq m_n$, $\sigma(i)\in A_n^*$. 

Let $h$ be a $\mathbb{Q}_{\alpha\beta}(g)$-generic over $L(H^G)$ such that $r_{(\xi)}\in h$. Then $h$ is $\mathbb{P}_{U_\beta}$-generic over $L(H^G)$. All dense subsets of $\mathbb{P}_{U_\beta}$ in $V$ are in $(H_{\kappa^+})^{V[G\downarrow\gamma]}$ for some $\gamma<\theta$, therefore $h$ is also $\mathbb{P}_{U_\beta}$-generic over $V$. Then, by Proposition \ref{prop:sequenceforforcing}, $\sigma_h$ is a branch of $T$ in $L(H^G)[h]$. But then, by absoluteness of well-foundedness, there is a $\sigma'\in L(H^G)$ branch of $T$, and by Proposition \ref{prop:sequenceforforcing} ${\cal F}_{\sigma'}$ is $\mathbb{P}_{U_\beta}$-generic over $V$. 

Consider $m_0\geq\lh(s)$ such that for all $i\geq m_0$, $\sigma'(i)\in A_0^*$. Then we have that $\sigma'\upharpoonright[m_0,\omega)\subseteq A\cap A_0^*$, so if we define $\overline{A}=(A\cap A_0^*)\setminus(\sigma'\upharpoonright m_0)$, then $p=(\sigma'\upharpoonright m_0,\overline{A})\leq r_{(\xi)}$ and every sequence compatible with $p$ is going to have a tail inside $A_0^*$. But this is still not enough, since it could be that $p$ is not $\alpha$-nice, i.e., $p\notin\mathbb{P}_{U_\beta}^\alpha$. By definition, though, $p\in{\cal F}_{\sigma'}$, that is $\mathbb{P}_{U_\beta}$-generic over $V$, and therefore by Lemma \ref{lem:nicedense}, there exists a $p_{(\xi)}\leq p$ such that $p_{(\xi)}\in\mathbb{P}_{U_\beta}^\alpha\cap{\cal F}_{\sigma'}$. Since $p_{(\xi)}\in{\cal F}_{\sigma'}$, $p_{(\xi)}\downarrow\alpha\in{\cal F}_{\sigma'\downarrow\alpha}={\cal F}_{\sigma_g}=g$, therefore $p_{(\xi)}\in\mathbb{Q}_{\alpha\beta}(g)$, and since $p_{(\xi)}\leq p\leq r_{(\xi)}$,  $p_{(\xi)}$ decides the value of $\tau(m)$ with $m\leq m_\emptyset$.

We continue by induction: for any $s\in\bigcup_{n\in\omega}\prod_{i<n}\eta_i$, say $\lh(s)=n$, suppose that $p_s\in\mathbb{Q}_{\alpha\beta}(g)$ is defined. Then there must exist an $m_s$ such that $|\Gamma_{m_s}(p_s)|\geq\eta_n$, so we can find in $\mathbb{Q}_{\alpha\beta}(g)$ at least $\eta_n$-many extensions of $p_s$ whose interpretations of $\tau$ differ on $m_s$. Choose $q_{s{}^\smallfrown(\xi)}$ to be the $\xi$-th of them, so that $q_{s{}^\smallfrown(\xi)}\Vdash\tau(\check{m}_s)=\check{\gamma}_{s{}^\smallfrown(\xi)}$ for some $\gamma_{s{}^\smallfrown(\xi)}$, so that if $\xi_1\neq\xi_2$ then $\gamma_{s{}^\smallfrown(\xi_1)}\neq\gamma_{s{}^\smallfrown(\xi_2)}$. Extend $q_{s{}^\smallfrown(\xi)}$ to $r_{s{}^\smallfrown(\xi)}$ so that $r_{s{}^\smallfrown(\xi)}$ can compute $\tau$ up until $m_s$ and $n+1$. Consider $T$ the tree of all generic sequences compatible with $r_{s{}^\smallfrown(\xi)}$: there is a branch in a generic extension, and therefore there is one in $L(H^G)$, and there must exist an $m_n$ such that its tail after $m_n$ is contained in $A^*_n$. From this, we can find $m_n^+$ and $p_{s{}^\smallfrown(\xi)}\in\mathbb{Q}_\alpha(g)$ such that its stem has length $m_n^+$, and its measure one part is inside $A^*_n$.

For each $x\in\prod_{i<\omega}\eta_i$, let $\sigma_x=\bigcup_{n<\omega}stem(p_{x\upharpoonright n})$. Note that for each $n$, $p_{x\upharpoonright n},q_{x\upharpoonright n},r_{x\upharpoonright n}\in{\cal F}_{\sigma_x}$. Let ${\cal F}^\alpha_{\sigma_x}={\cal F}_{\sigma_x}\cap\mathbb{P}_{U_\beta}^\alpha$. Because of our construction:
\begin{itemize}
 \item for each $x\in\prod_{i<\omega}\eta_i$, ${\cal F}^\alpha_{\sigma_x}$ is $\mathbb{Q}_{\alpha\beta}(g)$-generic over $V[g]$: for each $n\in\omega$, we have that the measure-one part of $p_{x\upharpoonright n}$ is inside $A^*_n$, therefore for each $i>m_n^+=\lh(p_{x\upharpoonright n})$ $\sigma_x(i)\in A^*_n$, so by Proposition \ref{prop:sequenceforforcing} ${\cal F}_{\sigma_x}$ is $\mathbb{P}_{U_\beta}$-generic over $V$ and ${\cal F}^\alpha_{\sigma_x}$ is $\mathbb{P}^\alpha_{U_\beta}$-generic over $V$. Since $p_{x\upharpoonright n}\in\mathbb{Q}_{\alpha\beta}(g)$, then $p_{x\upharpoonright n}\downarrow\alpha\in g$, so for each $n\in\omega$ $stem(p_{x\upharpoonright n})\downarrow\alpha\sqsubseteq\sigma_g$. But then $\sigma_x\downarrow\alpha=\sigma_g$, so ${\cal F}^\alpha_{\sigma_x}\downarrow\alpha={\cal F}_{\sigma_g}=g$, and ${\cal F}^\alpha_{\sigma_x}$ is $\mathbb{Q}_{\alpha\beta}(g)$-generic over $V[g]$;
 \item for each $x,y\in\prod_{i<\omega}\eta_i$, if $x\neq y$ then $\tau_{{\cal F}^\alpha_{\sigma_x}}\neq\tau_{{\cal F}^\alpha_{\sigma_y}}$: let $n$ be the maximal such that $x\upharpoonright n=y\upharpoonright n$. Then $q_{x\upharpoonright (n+1)}\Vdash\tau(\check{m}_{x\upharpoonright n})=\check{\gamma}_{x\upharpoonright (n+1)}$ and $q_{y\upharpoonright (n+1)}\Vdash\tau(\check{m}_{x\upharpoonright n})=\check{\gamma}_{y\upharpoonright (n+1)}$, with $\check{\gamma}_{x\upharpoonright (n+1)}\neq\check{\gamma}_{y\upharpoonright (n+1)}$;
 \item for each $x,y\in\prod_{i<\omega}\eta_i$, if $n>0$ is such that $x\upharpoonright n=y\upharpoonright n$, then there is an $m>n$ such that $\tau_{{\cal F}^\alpha_{\sigma_x}}\upharpoonright m=\tau_{{\cal F}^\alpha_{\sigma_y}}\upharpoonright m$: consider $r_{x\upharpoonright n}$; it is both in ${\cal F}^\alpha_{\sigma_x}$ and in ${\cal F}^\alpha_{\sigma_y}$, and forces a value for $\tau(i)$ for every $i<m_{x\upharpoonright (n-1)},n+1$, therefore for every $i<m_{x\upharpoonright (n-1)},n+1$ $\tau_{{\cal F}^\alpha_{\sigma_x}}(i)=\tau_{{\cal F}^\alpha_{\sigma_y}}(i)$.
\end{itemize}

The function $x\mapsto\tau_{{\cal F}^\alpha_{\sigma_x}}$ is therefore a continuous function from $C(\kappa)=\prod_{n\in\omega}\eta_n$ to $\{\tau_h : h$ is $\mathbb{Q}_\alpha(g)$-generic over $V[g]$, $h\in L(H^G)\}$, whose image is closed-in-$\pre{\omega}{\kappa}$, and therefore the theorem is proved.
\end{proof}

We are going to prove now that in $L(H^G)$ all subsets of $\pre{\omega}{\kappa}$ have the $\kappa$-PSP. The proof is similar to the original proof by Solovay of the consistency of ``All subsets of $\pre{\omega}{\omega}$ have the PSP'': Solovay split the forcing $Col(\omega,<\theta)$, with $\theta$ inaccessible, in $Col(\omega,<\alpha)\times Col(\alpha,<\beta)\times Col(\beta,<\theta)$, arguing that if $A$ is an uncountable and ordinal-definable set of reals in $V[G]$, a generic extension via $Col(\omega,<\theta)$, then by the homogeneity of $Col(\alpha,<\beta)$ there is a perfect set of $Col(\alpha,<\beta)$-generics over $V[G\upharpoonright\alpha]$, and by the homogeneity of $Col(\beta,<\theta)$, for all $Col(\beta,<\theta)$-generics $H$, each of the elements of the perfect set $x$ satisfy the formula that defines $A$ in $V[G\upharpoonright\alpha][x][H]$. But Solovay proved that there is an $H$ such that $V[G]=V[G\upharpoonright\alpha][x][H]$, so the perfect set is inside $A$.

In our case, we are going to split the forcing $\mathbb{P}_U$ in $\mathbb{P}_{U_\alpha}*\dot{\mathbb{Q}}_{\alpha\beta}*\dot{\mathbb{Q}}_\beta$ (it is not the same forcing, but it is equivalent), and we are going to use Theorem \ref{thm:perfect} in the first instance, and Corollary \ref{cor:homogeneity} in the second instance. It remains to show an analogue of the existence of an $H$ such that $V[G]=V[G\upharpoonright\alpha][x][H]$. This is what we are going to use: 

\begin{prop}{\cite[Proposition 3.1.7, Proposition 3.1.8]{Kafkoulis1994}}
 \label{prop:switchgenerics}
 Let $\alpha<\beta<\theta$. Then for each $h\in L(H^G)$ $\mathbb{Q}_{\alpha\beta}(G\downarrow\alpha)$-generic over $G[V\downarrow\alpha]$, there exists $G^*$ $\mathbb{P}_U$-generic over $V$ such that $\sigma_{G^*}\downarrow\beta=\sigma_h$ and $H^{G^*}=H^G$.\footnote{This result is announced at the beginning of Section 3.1, yet Proposition 3.1.8 is proved without taking into consideration $h$. The trick to reach the full result is in the proof of \cite[Theorem 3.2.1]{Kafkoulis1994}. In short, Kafkoulis defines a forcing notion $\mathbb{Q}$, and proves that each $\mathbb{Q}$-generic generates a $\mathbb{P}_U$-generic $G^*$ with the desired properties. If we take the $\mathbb{Q}$-generic that contains the condition $(\sigma_{h},0,\beta,{\cal P}_\kappa(\theta))\in\mathbb{Q}$, then $G^*$ is as we wanted.}
\end{prop}

So even if we have not proved that for each $G$ $\mathbb{P}_U$-generic, and for each $h$ $\mathbb{Q}_{\alpha\beta}(G\downarrow\alpha)$-generic, we can find a generic $H$ such that $V[G\downarrow\alpha][h][H]=V[G]$, anyway there is a $G^*$ such that 
\begin{equation*}
 L(V_{\kappa+1})^{V[G\downarrow\alpha][h][G^*]}=L(V_{\kappa+1})^{V[G^*]}=L(V_{\kappa+1})^{V[G]},
\end{equation*}
and this is enough.

\begin{teo} \label{thm:main}
 Each subset of $\pre{\omega}{\kappa}$ in $L(H^G)$ has the $\kappa$-PSP.
\end{teo}
\begin{proof}
  Let $A\subseteq\pre{\omega}{\kappa}$, $A\in L(H^G)$. Suppose that $(\neg(|A|\leq\kappa))^{L(H^G)}$. Since $A\in L(H^G)$, $A$ is definable in $L(H^G)$ with parameters in $\Ord\cup H^G\cup\{H^G\}$. Moreover, $L(H^G)$ is a definable class in $V[G]$ with parameter $H^G$, therefore there are a formula $\varphi$ and $x_0,\dots,x_n\in \Ord\cup H^G$ such that $x\in A$ iff $V[G]\vDash\varphi(x,x_0,\dots,x_n,H^G)$. Let $\alpha<\theta$ be such that $x_0,\dots,x_n\in V[G\downarrow\alpha]$, that exists by Proposition \ref{prop:basicH}. 

 The set $(\pre{\omega}{\kappa})^{V[G\downarrow\alpha]}$ has cardinality $\eta=2^\kappa$ in $V[G\downarrow\alpha]$. Since $\theta$ is inaccessible, there is a $\gamma<\theta$ such that $\gamma>\eta$, and $(|\gamma|=\kappa)^{V[G\downarrow\gamma]}$ by Theorem \ref{thm:basic}. So in $V[G\downarrow\gamma]$ there is a bijection between $(\pre{\omega}{\kappa})^{V[G\downarrow\alpha]}$ and $\kappa$, and this is codeable inside ${\cal P}(\kappa)\cap V[G\downarrow\gamma]$ \footnote{e.g. $\{(\alpha,n,\beta)|f(\alpha)(n)=\beta\}$, where $f$ is the bijection.}, therefore this is true also in $L(H^G)$. But then $(|(\pre{\omega}{\kappa})^{V[G\downarrow\alpha]}|=\kappa)^{L(H^G)}$. But since $(\neg(|A|\leq\kappa))^{L(H^G)}$, it must be that there is an $x\in A$ such that $x\notin (\pre{\omega}{\kappa})^{V[G\downarrow\alpha]}$. In particular $V[G]\vDash\varphi(x,x_0,\dots,x_n,H^G)$.

Let $\beta<\theta$ be such that $x\in V[G\downarrow\beta]$, that exists by Proposition \ref{prop:basicH}, and fix a $\mathbb{P}_{U}$-name $F$ for $H^G$ that is invariant under any $\pi_f$ as in Lemma \ref{lem:invariantname}. Then 
\begin{equation*}
 V[G\downarrow\beta][G\cap\mathbb{P}^\beta_U]\vDash\varphi(x,x_0,\dots,x_n,H^G),
\end{equation*}
so there is a $p\in\mathbb{Q}_\beta(G\downarrow\beta)\cap G$ such that $p\Vdash_{\mathbb{Q}_\beta(G\downarrow\beta)}\varphi(\check{x},\check{x}_0,\dots\check{x}_n,i(F^\beta)_{G\downarrow\beta})$, where $F^\beta$ is the $\mathbb{P}_U^\beta$-name equivalent to $F$ and $i$ is the dense embedding between $\mathbb{P}_U^\beta$ and $\mathbb{P}_{U_\beta}*\dot{\mathbb{Q}}_\beta$. Therefore, by Corollary \ref{cor:homogeneity}, $\Vdash_{\mathbb{Q}_\beta(G\downarrow\beta)}\varphi(\check{x},\check{x}_0,\dots\check{x}_n,i(F^\beta)_{G\downarrow\beta})$ is true in $V[G\downarrow\beta]$.

Since 
\begin{equation*}
 V[G\downarrow\beta]=V[G\downarrow\alpha][(G\downarrow\beta)\cap\mathbb{P}^\alpha_{U_\beta}],
\end{equation*}
there is a $\mathbb{Q}_{\alpha\beta}(G\downarrow\alpha)$-name $\tau$ for $x$, and there is a $p\in\mathbb{Q}_{\alpha\beta}(G\downarrow\alpha)$ that forces the above, so that for any $h$ $\mathbb{Q}_{\alpha\beta}(G\downarrow\alpha)$-generic over $V[G\downarrow\alpha]$ such that $p\in h$, 
\begin{equation*}
 V[G\downarrow\alpha][h]\vDash (\Vdash_{\mathbb{Q}_\beta(h)}\varphi((\tau_h)^{\vee},\check{x}_0,\dots\check{x}_n,i(F^\beta)_h)),
\end{equation*}
where $(\tau_h)^{\vee}$ is the canonical check-name of $\tau_h\in V[G\downarrow\beta]$ for the forcing $\mathbb{Q}_\beta(h)$.

Let $p_1\leq p$ be such that $p_1\Vdash\tau:\check{\omega}\to\check{\kappa}$. Then by Theorem \ref{thm:perfect} applied to $p_1$, we can find in $L(H^G)$ a set of $\mathbb{Q}_{\alpha\beta}(G\downarrow\alpha)$-generics $\{h_x : x\in\prod_{n\in\omega}\eta_n\}$, with $(\eta_n : n\in\omega)$ cofinal in $\kappa$, such that $\{\tau_{h_x} : x\in\prod_{n\in\omega}\eta_n\}$ is a $\kappa$-perfect subset of $\pre{\omega}{\kappa}$. Also, for each $x\in\prod_{n\in\omega}\eta_n$, $p\in h_x$, therefore 
\begin{equation*}
 V[G\downarrow\alpha][h_x]\vDash (\Vdash_{\mathbb{Q}_\beta(h_x)}\varphi((\tau_{h_x})^{\vee},\check{x}_0,\dots\check{x}_n,i(F^\beta)_{h_x})).
\end{equation*}

Fix $x\in\prod_{n\in\omega}\eta_n$. Then by Proposition \ref{prop:switchgenerics} there is a $G^*$ such that $\sigma_{G^*}\downarrow\beta=\sigma_{h_x}$ and $H^{G^*}=H^G$. Then $V[G^*]\vDash\varphi(\tau_{h_x},x_0,\dots,x_n,H^{G^*})$, and so $V[G^*]\vDash\varphi(\tau_{h_x},x_0,\dots,x_n,H^G)$. But $\varphi(\tau_{h_x},x_0,\dots,x_n,H^G)$ says: $\tau_{h_x}$ is an element of $\pre{\omega}{\kappa}$ in $L(H^G)$ that satisfies in $L(H^G)$ the formula that witnesses the membership in $A$ using parameters $x_0,\dots,x_n$, and this is absolute among models that contain $H^G$ and $x_0,\dots,x_n$, therefore $V[G]\vDash\varphi(\tau_{h_x},x_0,\dots,x_n,H^G)$, and therefore $\tau_{h_x}\in A$ for all $x\in\prod_{n\in\omega}\eta_n$. The map $x\mapsto \tau_{h_x}$ then witnesses that $A$ satisfies the $\kappa$-PSP.
\end{proof}

\begin{cor}
 \begin{center}
   Con($\ZFC+\exists\kappa,\theta (\kappa<\theta$, $\theta$ is inaccessible and $\kappa$ is a $\theta$-supercompact cardinal))
	
	 $\rightarrow$
	
	 Con($\ZF+\exists\kappa>\omega (\DC_\kappa+$ there is a $\kappa$-Polish space $X$ such that all the subsets of $X$ have the $\kappa$-PSP)) (P).
 \end{center}
\end{cor}

Note that if $V$ is a model of $\ZFC$, then $L(V_{\lambda+1})$ is a model of $\ZF+\DC_\lambda$ (see e.g. \cite[Lemma 4.10]{Dimonte2018}), therefore the previously known upperbound for the consistency strength of (P) was I0.

In fact, we can substitute ``there is a $\kappa$-Polish space'' in (P) with ``For all $\kappa$-Polish spaces'' thanks to the following Corollary:

\begin{cor}
 Let $X$ be a $\kappa$-Polish space in $L(H^G)$. Then all its subsets in $L(H^G)$ have the $\kappa$-PSP.
\end{cor}
\begin{proof}
 Let $X$ be a $\kappa$-Polish space in $L(H^G)$ and $A\subseteq X$. Since $X$ is $\kappa$-Polish, there is a closed set $F\subseteq\pre{\omega}{\kappa}$ and a continuous bijection $f:F\to X$. Consider $B=f^{-1}(A)$. Then $B\subseteq\pre{\omega}{\kappa}$ and $B\in L(H^G)$, therefore by Theorem \ref{thm:main} $B$ has the $\kappa$-PSP. If $|B|\leq\kappa$, then $|A|\leq\kappa$, and if $\pre{\kappa}{2}$ embeds into $B$, then $\pre{\kappa}{2}$ embeds into $A$.
\end{proof}

\end{document}